%Format: plain

%% arXiv:0812.3258

\input amstex
\documentstyle{amsppt}

\input label.def
\input degt.def
%\input debug.def
%\PrintLabels

%\copycounter\thm\subsection
%\copycounter\equation\subsection

\input epsf
\def\picture#1{\epsffile{#1-bb.eps}}

\def\ie.{\emph{i.e\.}}
\def\eg.{\emph{e.g\.}}
\def\cf.{\emph{cf\.}}

%% From the previous paper %%

%%%%%%%%%%%%%%%%%%%%%%%%%%%%%

{\catcode`\@11
\gdef\proclaimfont@{\sl}}

\newtheorem{Remark}\Remark\thm\endAmSdef
\conjecture\thm\endproclaim
%\newhead\subsubsection\subsubsection\endsubhead
%\def\paragraph{\subsubsection{}}
\def\paragraph{\subsection{}}
\copycounter\thm\subsection
\copycounter\equation\subsection

\loadbold
\def\bA{\bold A}

\def\bE{\bold E}

\def\tA#1{\tilde\bA#1}

\def\tE#1{\tilde\bE#1}

\def\bary{\bar y}
\def\barz{\bar z}

\let\splus\oplus

\let\into\hookrightarrow
\let\onto\twoheadrightarrow

\def\CG#1{\Z_{#1}}
\def\BG#1{\Bbb B_{#1}}
\def\DG#1{\Bbb D_{#1}}

\let\Gs\sigma

\def\CC{\Cal C}

\def\tX{\tilde X}
\def\tF{\tilde F}
\def\tpi{\tilde\pi}
\def\tm{\tilde m}

\def\B{\bar B}

\def\Cp#1{\Bbb P^{#1}}
\def\term#1-{$\DG{#1}$-}

\let\Ga\alpha

\let\Gg\gamma

\let\Gs\sigma

\def\1{^{-1}}

\def\ls|#1|{\mathopen|#1\mathclose|}
\let\<\langle
\let\>\rangle

\def\Aut{\operatorname{Aut}}
\def\Sk#1{\operatorname{Sk}#1}

\def\0{\text{--\,}}

\def\fragment(#1)#2{\ref{fig.e8}(#1)--$#2$}
\def\bifrag(#1)#2{\fragment(#1){#2,\bar#2\!\!\!}}
\def\rfragment(#1)#2{\ref{fig.e8-r}(#1)--$#2$}
\def\rbifrag(#1)#2{\rfragment(#1){#2,\bar#2\!\!\!}}

\def\inserthyphen{\ifcat\next a-\fi\ignorespaces}
\def\pblack-{$\bullet$\futurelet\next\inserthyphen}
\def\pwhite-{$\circ$\futurelet\next\inserthyphen}
\def\pcross-{$\vcenter{\hbox{$\scriptstyle\times$}}$\futurelet\next\inserthyphen}
\def\black{\protect\pblack}

\def\beginGAP{\bgroup
 \catcode`\^=12\catcode`\#=12\catcode`\_=12
 \obeylines\obeyspaces\eightpoint\tt}
\let\endGAP\egroup

\topmatter

\author
Alex Degtyarev
\endauthor

\title
Classical Zariski pairs
\endtitle

\address
Department of Mathematics,
Bilkent University,
06800 Ankara, Turkey
\endaddress

\email
degt\@fen.bilkent.edu.tr
\endemail

\abstract
We compute the fundamental groups of
%the complements of
all irreducible plane sextics
constituting classical Zariski pairs
\endabstract



\endtopmatter

\document

\section{Introduction}

A \emph{classical Zariski pair} is a pair of irreducible plane
sextics that share the same combinatorial type of singularities
but differ by the Alexander polynomial~\cite{Libgober1}.
The first example of such a pair was constructed by
O.~Zariski~\cite{Zariski.group}. Then, it was shown in~\cite{poly}
that the curves constituting a classical Zariski pair have simple
singularities only and, within each pair, the Alexander polynomial
of one of the curves is $t^2-t+1$, whereas the polynomial of the
other curve is trivial. The former curve is called
\emph{abundant}, and the latter \emph{non-abundant}. The abundant
curve is necessarily of \emph{torus type}, \ie., its equation can
be represented in the form $f_2^3+f_3^2=0$, where $f_2$ and~$f_3$
are homogeneous polynomials of degree~$2$ and~$3$, respectively.

A complete classification of classical Zariski pairs up to
equisingular deformation was recently obtained by
A.~\"Ozg\"uner~\cite{Aysegul}. Altogether,
there are $51$
%classical Zariski pairs of irreducible sextics (
pairs, one of them
being in fact a triple: the non-abundant curves with the set of
singularities $\bE_6\oplus\bA_{11}\oplus\bA_1$ form two distinct
complex conjugate deformation families. The purpose of this note is
to compute the fundamental groups of (the complements of) the
curves constituting classical Zariski pairs. We prove the
following theorem.

\theorem\label{th.main}
Within each classical Zariski pair, the fundamental group of the
abundant \rom(respectively, non-abundant\rom) curve is
$\BG3/(\Gs_1\Gs_2)^3$ \rom(respectively, $\CG6$\rom).
\endtheorem

This theorem is proved in Section~\ref{S.proof}, using the list
of~\cite{Aysegul} and a case by case analysis. In fact, most
groups are already known, see~\cite{Artal.Trends},
\cite{degt.2a8}, \cite{dessin-e7}, \cite{dessin-e8},
and~\cite{Oka.conjecture},
and the
%only missing case is the set of singularities
%$\bA_{14}\splus\bA_2\splus2\bA_1$. We obtain the missing curves by
few missing curves can be obtained by
perturbing the set of singularities $\bA_{17}\splus2\bA_1$. The
construction and the computation of the fundamental group are
found in Sections~\ref{S.non-torus} (the non-abundant curves)
and~\ref{S.torus} (the abundant curves).

\section{The curve not of torus type\label{S.non-torus}}

\paragraph\label{s.non-torus}
Up to projective transformation, there is a unique curve
$C\subset\Cp2$ with the set of singularities
$\bA_{17}\splus2\bA_1$ and not of torus type, see~\cite{Shimada};
its transcendental lattice is
$\Bigl[\smallmatrix4&2\\2&10\endsmallmatrix\Bigr]$.
After nine blow-ups, the curve transforms to the union of
two of the three type~$\tA{_0^*}$ fibers in a
Jacobian rational elliptic
surface with
%a type~$\tA{_8}$ fiber.
the combinatorial type of singular fibers
$\tA{_8}\splus3\tA{_0^*}$.
For the equation, consider
the pencil of cubics given by
$$
f_b(x,y):=b(-x^2-xy^2+y)+(x^3-xy+y^3)=0,\qquad b\in\Cp1,
$$
and take two fibers corresponding to two distinct roots of
$b^3=1/27$. (All three roots give rise to
nodal cubics, which are the three
type~$\tA{_0^*}$ fibers in the elliptic pencil above. The curve
corresponding to $b=\epsilon/3$, $\epsilon^3=1$, has a node at
$x=(2/5)\epsilon\1$, $y=(1/5)\epsilon$. The type~$\tA{_8}$ fiber
blows down to the nodal cubic $\{f_0=0\}$.)

\lemma\label{group.non-torus}
For the curve~$C$ as in~\ref{s.non-torus}, one has
$$
\pi_1(\Cp2\sminus C)=\<p,\Gg_+\,|\,
 \text{$p^9=1$, $\Gg_+\1p\Gg_+=p^4$}\>.
$$
\endlemma

\proof
Consider the trigonal curve $\B\subset\Sigma_2$ with a type~$\bA_8$
singular point. Its skeleton~$\Sk$,
see~\cite{degt.kplets}, is shown in
Figure~\ref{fig.a8}.

\midinsert
\centerline{\picture{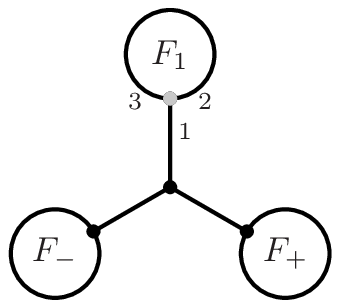}}
\figure
The skeleton~$\Sk$ of~$\B$
\endfigure\label{fig.a8}
\endinsert

Let $F_1$, $F_\pm$ be the type~$\tA{_0^*}$
singular fibers of~$\B$, and let $F_\infty$ be the type~$\tA{_8}$
fiber.
(Recall that $F_1$, $F_\pm$ are located inside the small loops in
Figure~\ref{fig.a8}, whereas $F_\infty$ is inside the outer
region.)
Consider the
minimal resolution of the
double covering $\tX\to\Sigma_2$ ramified
at~$\B$ and the exceptional section $E\subset\Sigma_2$,
and denote by tildes the
pull-backs of the fibers in~$\tX$.

Consider the nonsingular fiber~$F$ over the %central
\black-vertex~$v$ of~$\Sk$ next to~$F_1$ (shown in grey in
Figure~\ref{fig.a8}),
denote
$\pi_F:=\pi_1(F\sminus(\B\cup E))$,
and pick a canonical basis
$\{\Ga_1,\Ga_2,\Ga_3\}$ for~$\pi_F$
defined by the marking of~$\Sk$ at~$v$ shown in
%the figure,
Figure~\ref{fig.a8},
see~\cite{degt.kplets}.
Then the
fundamental
group
$\tpi_F:=\pi_1(\tF\sminus E)$ of the punctured torus
$\tF\sminus E$
is obtained from~$\pi_F$ by adding
the relations $\Ga_1^2=\Ga_2^2=\Ga_3^2=1$ and passing to the
kernel of the homomorphism $\pi_F\to\CG2$,
$\Ga_1,\Ga_2,\Ga_3\mapsto1$. Hence, $\tpi_F$ is the free group
generated by
$$
p:=\Ga_1\Ga_2=(\Ga_2\Ga_1)\1\quad\text{and}\quad
q:=(\Ga_3\Ga_2)=(\Ga_2\Ga_3)\1.
$$

Start with the group
$$
G_1=\pi_1(\tX\sminus(E\cup\tF_+\cup\tF_-\cup\tF_\infty))
$$
and compute it
applying Zariski--van Kampen's approach~\cite{vanKampen} to the
elliptic pencil on~$\tX$.
Let $\Gg_1$, $\Gg_\pm$ be the
generators of the free group
$$
\pi_1(\Cp1\sminus(F_1\cup F_+\cup F_-\cup F_\infty), F)
$$
represented by the shortest loops in~$\Sk$ starting at~$v$ and
circumventing the corresponding fibers in the counterclockwise
direction.
(We identify fibers of the ruling and their projections to the
base.)
Using a proper section, see~\cite{degt.kplets},
one can lift these generators to
$\Sigma_2\sminus(\B\cup E)$ and to $\tX\sminus E$. Using the same
proper section, define the braid monodromies
$m_1,m_\pm\in\Aut\pi_F$ and their lifts
$\tm_1,\tm_\pm\in\Aut\tpi_F$. In this notation, the group~$G_1$
has the following presentation, \cf.~\cite{vanKampen}:
$$
G_1=\bigl<p,q,\Gg_+,\Gg_-\bigm|
 \text{$p=\tm_1(p)$, $q=\tm_1(q)$,
 $\Gg_\pm\1p\Gg_\pm=\tm_\pm(p)$,
 $\Gg_\pm\1q\Gg_\pm=\tm_\pm(q)$}\bigr>.
$$
The braid monodromy is computed as explained
in~\cite{degt.kplets};
for~$\B$ it is
$$
m_1=\Gs_2,\quad
m_+=\Gs_1^{-3}\Gs_2\Gs_1^3,\quad
m_-=\Gs_1\1\Gs_2^2\Gs_1\Gs_2^{-2}\Gs_1,
$$
where $\Gs_1$, $\Gs_2$ are the Artin generators of~$\BG3$
(we assume that the braid group~$\BG3$ acts on~$\pi_F$ from the
left), and in terms of~$p$ and~$q$ it takes the form
$$
\gather
\tm_1\:p\mapsto pq,\quad q\mapsto q;\\
%%p := [1, 2, 3, -2, -1, -2, -1, 2, 1, 2]
%%q := [-2, -1, -2, 1, 2, 1, 2, 1, 2, 3, -2, -1, -2, -1, 2, 1, 2, 1, 2, -3, -2, -1]
\tm_+\: p\mapsto pqp^3,\quad q\mapsto p^{-4}q\1p^{-4}q\1p\1;\\
%%p := [1, 2, 3, -2, 1, 2, -3, -2, -1, -2, 1, 2, -3, -2, -1, 2, 1, 2, 3, -2, -1, 2]
%%q := [-2, 1, 2, -3, -2, -1, -2, 1, 2, 3, -2, -1, 2, 1, 2, 3, -2, -1, 2, 1, 2, 3, -2, 1, 2, -3, -2, -1]
\tm_-\: p\mapsto(pq)^2(p^2q)^2p,\quad
q\mapsto p\1q\1(p^{-2}q\1)^3p\1q\1p\1.
\endgather
$$
%Now, the
The
very first relation $p=pq$ implies $q=1$. Hence
also $\tm_\pm(q)=1$ and $p^9=1$. Thus, one has
$$
G_1=\bigl<p,\Gg_+,\Gg_-\bigm|
 \text{$p^9=1$, $\Gg_+\1p\Gg_+=p^4$,
 $\Gg_-\1p\Gg_-=p^7$}\bigr>.
\eqtag\label{eq.G1}
$$

In order to pass to the group $\pi_1(\Cp2\sminus B)$,
we need to patch
back in one of the nine irreducible
components of the type~$\tA{_8}$ fiber~$F_\infty$.
(The component to be patched in
is the proper transform of the nodal curve
$\{f_0(x,y)=0\}$.) This operation adds to~\eqref{eq.G1} an
additional relation $[\partial\tilde\Gamma]=1$, where
$\tilde\Gamma$ is a small holomorphic disk in~$\tX$ transversal to
the component in question.
Using a proper section again, one can see that in~$G_1$ there is
a relation
$[\partial\tilde\Gamma]\1p^?=\Gg_-\Gg_+$, where $p^?$ is merely
an element of the group~$\tpi_F$
of the fiber (modulo the relations in~$G_1$),
which we do not bother to
compute. Adding the extra relation
$[\partial\tilde\Gamma]=1$ to~\eqref{eq.G1} and
eliminating~$\Gg_-$, one arrives at the presentation
announced in the
statement.
(Note that $7=4\1\bmod9$, hence the order of~$p$
remains~$9$.)
\endproof

\corollary
The commutant of the group $\pi_1(\Cp2\sminus C)$ as in
Lemma~\ref{group.non-torus} is a
central subgroup of order~$3$.
\endcorollary

\proof
The commutant is normally generated by the commutator
$p\1\Gg_+\1p\Gg_+=p^3$; it is a central element of order~$3$.
\endproof

\corollary\label{cor.Z6}
For any irreducible perturbation $C'$ of the curve~$C$ as
in~\ref{s.non-torus}, one has $\pi_1(\Cp2\sminus C')=\CG6$.
\endcorollary

\proof
Any central extension $\{1\}\to\CG3\to G\to\CG6\to\{1\}$
would be abelian.
\endproof

\section{The curve of torus type\label{S.torus}}

\paragraph\label{s.torus}
Up to projective transformation, there is a unique torus type curve
$C\subset\Cp2$ with the set of singularities
$\bA_{17}\splus2\bA_1$, see~\cite{Shimada};
its transcendental lattice is
$\Bigl[\smallmatrix2&0\\0&2\endsmallmatrix\Bigr]$.
Similar to~\ref{s.non-torus}, this curve blows up to the union of
the two type $\tA{_0^*}$ fibers in a Jacobian rational elliptic
surface with the combinatorial type of singular fibers
$\tE{_8}\splus2\tA{_0^*}$.
The curve can be given by the equation
$$
f(x,y):=(y^3+y^2+x^2)\Bigl(y^3+y^2+x^2-\frac4{27}\Bigr)=0,
$$
and its torus structure is
$$
f(x,y)=
 \Bigl(y^3+y^2+x^2-\frac2{27}\Bigr)^2+\Bigl(\frac{\root3\of4}{9}\Bigr)^3.
$$

\lemma\label{group.torus}
Let $C$ be a curve as in~\ref{s.torus}, and let~$U$ be a Milnor
ball about the type~$\bA_{17}$ singular point of~$C$. Then the
homomorphism $\pi_1(U\sminus C)\to\pi_1(\Cp2\sminus C)$
induced by the inclusion $U\into\Cp2$
is surjective.
\endlemma

\proof
In the coordinates $\bary=y/x$, $\barz=1/x$, the curve is given by
the equation
$$
(\bary^3+\bary^2\barz+\barz)
 \Bigl(\bary^3+\bary^2\barz+\barz-\frac4{27}\barz^3\Bigr)=0,
$$
the type~$\bA_{17}$ singular point is at the origin, and each
component is inflection tangent to the line $\{\barz=0\}$ at this
point.
To compute the group, apply Zariski--van Kampen
theorem~\cite{vanKampen} to the vertical pencil
$\{\barz=\const\}$, choosing for the reference
a generic fiber $F=\{\barz=\epsilon\}$ close to the
origin. On the one hand, one has an epimorphism
$\pi_1(F\sminus C)\onto\pi_1(\Cp2\sminus C)$. On the
other hand, the intersection $C\cap\{\barz=0\}$ consists of a
single $6$-fold point; hence, if $\epsilon$ is small enough, all
six points of the intersection $C\cap F$ belong to~$U$ and the
generators of $\pi_1(F\sminus C)$ can be chosen inside~$U$.
\endproof

\corollary\label{cor.B3}
Let~$C'$ be a perturbation of the curve~$C$ as in~\ref{s.torus}
with the set of singularities $\bA_{14}\splus\bA_2\splus2\bA_1$.
Then $\pi_1(\Cp2\sminus C')=\BG3/(\Gs_1\Gs_2)^3$.
\endcorollary

\proof
%Any perturbation as in the statement is necessarily irreducible
%and of torus type.
Let~$U$ be as in Lemma~\ref{group.torus}. Then
$\pi_1(U\sminus C')=\BG3$ and, due to the lemma, there is an
epimorphism $\BG3\onto\pi_1(\Cp2\sminus C')$. Since $C'$ is
necessarily irreducible and of torus type
(so that the
abelianization of $\pi_1(\Cp2\sminus C')$ is~$\CG6$ and
$\pi_1(\Cp2\sminus C')$ factors to $\BG3/(\Gs_1\Gs_2)^3$),
the latter epimorphism factors through an isomorphism
$\BG3/(\Gs_1\Gs_2)^3\cong\pi_1(\Cp2\sminus C')$.
\endproof

\Remark
The other irreducible perturbations of~$C$ that are of torus type
are considered elsewhere, see~\cite{degt.2a8}. Their groups are
also $\BG3/(\Gs_1\Gs_2)^3$.
\endRemark

\section{Proof of Theorem~\ref{th.main}\label{S.proof}}

\paragraph\label{pf.torus}
The groups of all but one sextics of torus type occurring in
classical Zariski pairs are known, see~\cite{degt.2a8} for a `map'
and further references; all groups are $\BG3/(\Gs_1\Gs_2)^3$.
The only missing curve has the set of singularities
$\bA_{14}\splus\bA_2\splus2\bA_1$.
Such a curve can be obtained by a
perturbation from a reducible sextic of torus type with the set of
singularities $\bA_{17}\splus2\bA_1$ (see Proposition~5.1.1
in~\cite{degt.8a2}), and its group is given by
Corollary~\ref{cor.B3}.

\paragraph\label{pf.non-torus}
The fundamental groups of most non-abundant sextics appearing in
classical Zariski pairs are computed in~\cite{degt.2a8},
\cite{dessin-e7},~\cite{dessin-e8}, with a
considerable contribution
from~\cite{Oka.conjecture}. According to~\cite{dessin-e8}, unknown
are the groups of the curves with the sets of singularities
$$
\bA_{17}\splus\bA_1,\quad
\bA_{14}\splus\bA_2\splus2\bA_1,\quad %(torus type is unknown either),
2\bA_8\splus2\bA_1,\quad
2\bA_8\splus\bA_1.
$$
The first curve can be obtained by a perturbation from a sextic
with a single type $\bA_{19}$ singular point. According
to~\cite{Artal.Trends}, its group is abelian. The three other
curves are perturbations of the curve~$C$ constructed
in~\ref{s.non-torus}, and their groups are abelian due to
Corollary~\ref{cor.Z6}.
(Note that the perturbations
%used
exist due to Proposition 5.1.1
in~\cite{degt.8a2}, and the resulting curves are unique up to
equisingular deformation due to~\cite{Aysegul}.)
\qed

\Remark
A curve~$C$ as in~\ref{s.non-torus} can also be perturbed to a
sextic with the set of
singularities $\bA_{17}\splus\bA_1$, but the result is reducible.
\endRemark

\refstyle{C}

\enddocument